%% file: main.tex
\documentclass[12pt,a4paper]{article}
\usepackage[biblatex]{anziamjedraft}
\usepackage{amsfonts}
\usepackage{todonotes}
\bibliography{references}
\usepackage{subcaption}
\usepackage{tikz}
\usepackage{pgfplots}
\usepackage{soul}
\usepackage{xcolor}

\pgfplotsset{compat=newest}
\usepgfplotslibrary{groupplots}
\usepgfplotslibrary{dateplot}

\title{On well-posed boundary conditions and energy stable finite volume method for the linear shallow water wave equation}

\author{Rudi Prihandoko}
\address{Mathematical Science Institute, Australian National University, 
Australian Capital Territory~2600, \textsc{Australia}.}
\mailto{u7107271@anu.edu.au}
\myorcid{'0000-0001-6376-7952'}
\author{Kenneth Duru}
\address{Mathematical Science Institute, Australian National University, 
Australian Capital Territory~2600, \textsc{Australia}.}
\author{Stephen Roberts}
\address{Mathematical Science Institute, Australian National University, 
Australian Capital Territory~2600, \textsc{Australia}.}
\author{Christopher Zoppou}
\address{Mathematical Science Institute, Australian National University, 
Australian Capital Territory~2600, \textsc{Australia}.}
\begin{document}

\maketitle

\begin{abstract}
	We derive and analyse well-posed boundary conditions for the linear shallow water wave equation. The analysis is based on the energy method and it identifies the number, location and form of the boundary conditions so that the initial boundary value problem is well-posed. A finite volume method is developed based on the summation-by-parts framework with the boundary conditions implemented weakly using penalties. Stability is proven by deriving a discrete energy estimate analogous to the continuous estimate.   The  continuous and discrete analysis covers all flow regimes. 
    Numerical experiments are presented verifying the analysis.
\end{abstract}

\tableofcontents

\section{Introduction}
\label{sec:intro}
Numerical models that solve the shallow water wave equations (SWWE) have become a common tool for modeling environmental problems. This system of nonlinear hyperbolic partial differential equations (PDE) represent the conservation of mass and momentum of unsteady free surface flow subject to gravitational forces. The SWWE  assume that the fluid is inviscid, incompressible and the wavelength of the wave is much greater than its height. Typically these waves are associated with flows caused for example by tsunamis, storm surges and floods in riverine systems. The SWWE are also a fundamental component for predicting a range of aquatic processes, including sediment transport and the transport of pollutants. All these processes can have a significant impact on the environment, vulnerable communities and infrastructure. Therefore making accurate predictions using the SWWE crucial for urban, rural and environmental planners. 

For practical problems, the SWWE has been solved numerically using finite difference methods \cite{mahmood1975unsteady}, finite volume methods \cite{zoppou2003explicit}, discontinuous Galerkin method \cite{WINTERS2015357} and the method of characteristics \cite{cunge1980practical}.  Although, the shallow water wave equations are in common use, a rigorous theoretical investigation of boundary conditions necessary for their solution is still an area of active research \cite{GHADER20141}.

In this paper, we investigate well-posed boundary conditions for the linearized SWWE using the energy method \cite{gustafsson2007high,gustafsson1995time} and develop provably stable numerical method for the model. Following Ghader and Nordstr\"{om} \cite{GHADER20141}, our analysis identifies the type, location and number of boundary conditions that are required to yield a well-posed initial boundary value problem (IBVP). More importantly, we formulate the boundary conditions so that they can be readily implemented in a stable manner for numerical approximations that obey the summation-by-parts (SBP)  principle \cite{kreiss1974finite}. 
We demonstrate this by deriving a stable finite volume method using the SBP framework and impose the boundary conditions weakly using the Simultaneous Approximation Term (SAT) method \cite{carpenter1994time}. This SBP-SAT approach enables us to prove that the numerical scheme satisfies the discrete counterparts of energy estimates required for well-posedness of the IBVP, resulting in a provably stable and conservative numerical scheme.

The continuous and discrete analysis covers all flow regimes, namely sub-critical, critical and super-critical flows.  Numerical experiments are performed to verify the theoretical analysis of the continuous and discrete models. 
\section{Continuous analysis}
\label{sec:con-analysis}
    The one dimensional SWWE are
    \begin{align}
        \label{eq:swwe}
        &\frac{\partial h}{\partial t} 
            + \frac{\partial (uh)}{\partial x} = 0, \quad \frac{\partial (uh)}{\partial t} 
            + \frac{\partial (u^2 h + \frac12 gh^2) }{\partial x} = 0,
    \end{align}
where $x\in \mathbb{R}$ is spatial variable, $t\ge 0$ is time, $h(x,t)>0$ and $u(x,t)$ are the water depth and the depth averaged fluid velocity respectively,  $g>0$ is the gravitational acceleration. 
    
To make our analysis tractable we linearise the SWWE by substituting $h = H + \widetilde{h}$ and $u = U+\widetilde{u}$ into \eqref{eq:swwe}, where $\widetilde{h}$ and $\widetilde{u}$ denote perturbations of the constant water depth $H>0$ and fluid velocity $U$ respectively.

After simplifying, the linearised SWWEs are
    \begin{align}
     \label{eq:lswwe}
        &\frac{\partial h}{\partial t} + U \frac{\partial h}{\partial x} 
            + H \frac{\partial u}{\partial x} = 0, \quad 
        \frac{\partial u}{\partial t} + g \frac{\partial h}{\partial x} 
            + U \frac{\partial u}{\partial x} = 0,
    \end{align}
where we have dropped the tilde on the perturbed variables.
    
Introducing the unknown vector field $\textbf{q} = \begin{bmatrix}h, &u\end{bmatrix}^\top$, the linear equation \eqref{eq:lswwe} can be rewritten in a more compact form as
    \begin{align}\label{eq:lswwe_compact}
        \frac{\partial \textbf{q}}{\partial t} = D\textbf{q}, \quad {D} = - {M}\frac{\partial}{\partial x}, \quad {M}=\begin{bmatrix}U&H\\g&U\end{bmatrix}.
    \end{align}
%
We will consider \eqref{eq:lswwe_compact} in a bounded domain and augment it with initial and boundary conditions.
Let our domain be $\Omega=[0,1]$ and $\Gamma =\{0,1\}$ be the boundary points.
    We consider the IBVP  
    \begin{subequations}\label{eq:ibvp}
        \begin{align}
            \frac{\partial \textbf{q}}{\partial t} &= D\textbf{q} \label{eq:ibvp_pde}, \ x\in \Omega, \ t\ge 0,\\
            \textbf{q}(x,0) &= \textbf{f}(x), \ x\in \Omega,\\
            \mathcal{B}\textbf{q} &=  \textbf{b}(t),  \ x\in \Gamma, \ t\ge 0,\label{eq:ibvp_boundary}
        \end{align}
    \end{subequations}
    where $\mathcal{B}$ is a linear boundary operator, $\textbf{b}$ is the boundary data and $\textbf{f} \in L^2(\Omega)$ is the initial condition. 
    One objective of this study is to investigate the choice of $\mathcal{B}$ which ensure that the IBVP \eqref{eq:ibvp} is well-posed. To simplify  the coming analysis, we will consider zero boundary data $\textbf{b}=0$, but the results can be extended to nontrivial boundary data $\textbf{b} \ne 0$. Furthermore, numerical experiments performed later in this paper confirm that the analysis is valid for nonzero boundary data.

        Let $\textbf {p}$ and $\textbf {q}$ be real valued functions, and define the weighted scalar product and the norm 
        \begin{align}
            \label{eq:L2_norm}
            \left(\textbf{p},\textbf{q}\right)_{W} = \int_\Omega \textbf{p}^\top {W}\textbf{q} \,\text{d}x, \qquad 
            \Vert \textbf{q} \Vert^2_{W} = (\textbf{q},\textbf{q})_{W},
        \end{align}
where $W=W^\top$  and $\textbf{q}^\top W\textbf{q} >0$
for all non-zero $\textbf{q} \in \mathbb{R}^2$. If $W=I$ we get the standard $L_2$ scalar product, and we omit the subscript $W$.


    \begin{definition}
        The  IBVP \eqref{eq:ibvp} is  well-posed if a unique solution $\textbf{q}$  satisfies
        \begin{align}
            \Vert \textbf{q}(\cdot,t) \Vert_{W} \leq \kappa e^{\nu t} \Vert \textbf{f} \Vert_{W}, \quad \Vert \textbf{f} \Vert_{W} < \infty,
        \end{align}
        for  some constants $\kappa>0$ and $\nu \in \mathbb{R}$ independent of $\textbf{f}$.
    \end{definition}
    The well-posedness of the IBVP \eqref{eq:ibvp} can be related to the boundedness of the differential operator ${D}$. We introduce the function space
     \begin{align}\label{eq:function_space}
        \mathbb{V} =\{\textbf{p}| \quad \textbf{p}(x)\in \mathbb{R}^2, \quad \Vert \textbf{p} \Vert_{W} < \infty, \quad 0\le x\le 1, \quad \{\mathcal{B}\textbf{p}=0,  \ x\in \Gamma \} \}.
       \end{align} 
    The following two definitions are useful.
    %
    %
    \begin{definition}\label{def:semi-boudedness}
        The  operator ${D}$
        is said to be \textbf{semi-bounded}  in 
        the function space $\mathbb{V}$ if it satisfies 
        \begin{align}
            (\textbf{q},{D}\textbf{q})_{W} 
                \leq \nu \Vert \textbf{q} \Vert^2_{W}, \quad \nu \in \mathbb{R}.
        \end{align}
         \end{definition}
          \begin{definition}
         The differential operator ${D}$ is \textbf{maximally semi-bounded} if it is semi-bounded in 
        the function space $\mathbb{V}$
        but not semi-bounded in any space 
        with fewer boundary conditions. 
    \end{definition}
    It is well-known that the {maximally semi-boundedness} of differential operator ${D}$ is a necessary and sufficient condition for the well-posedness of the IBVP \eqref{eq:ibvp} \cite{gustafsson1995time}. 
Thus to ensure that the IBVP \eqref{eq:ibvp} is well-posed, we need: a)  the differential operator ${D}$ to be semi-bounded and; b) the minimal number of boundary conditions  such that ${D}$ is maximally semi-bounded. 

To begin, we will show that the differential operator ${D}$ is semi-bounded in a certain weighted $L_2$ scalar product.
    \begin{lemma}\label{Lem:Semi-boundedness}
Consider the differential operator $D$ with the constant coefficients matrix $M$ given in \eqref{eq:lswwe_compact} and the  weighted $L_2$ scalar product defined in \eqref{eq:L2_norm}, where $W=W^\top$  and $\textbf{q}^\top W\textbf{q} >0$
for all non-zero $\textbf{q} \in \mathbb{R}^2$. If the matrix product $\widetilde{M} = WM$ is symmetric, $\widetilde{M} = \widetilde{M}^T$, and $ \left.\left(\textbf{q}^\top \widetilde{M} \textbf{q}\right)\right|_0^1 \ge 0$,   then $D$ is semi-bounded.
    \end{lemma}
%

    \begin{proof}
    We consider $\left(\textbf{q},\textbf{D}\textbf{q}\right)_{W}$ and use integration-by-parts, we have
    {
    \begin{align*}
        \left(\textbf{q},\textbf{D}\textbf{q}\right)_{W} 
        = -\int_\Omega \textbf{q}^\top \widetilde{M}\frac{\partial \textbf{q}}{\partial x} \,\text{d}x = -\frac{1}{2}\int_\Omega \frac{\partial }{\partial x} \left(\textbf{q}^\top \widetilde{M} \textbf{q}\right) \,\text{d}x = -\left.\frac{1}{2}\left(\textbf{q}^\top \widetilde{M} \textbf{q}\right)\right|_0^1.
    \end{align*}
    }
    Thus if  the boundary term $ \left.\left(\textbf{q}^\top \widetilde{M} \textbf{q}\right)\right|_0^1 \ge 0$ then $\left(\textbf{q},\textbf{D}\textbf{q}\right)_\textbf{W} \le  0$. In particular the lower bound $\left(\textbf{q},\textbf{D}\textbf{q}\right)_\textbf{W} =  0$ satisfies Definition \ref{def:semi-boudedness} with $\nu =0$.
    \end{proof}
    
    The next step will be to derive boundary operators $\{\mathcal{B}\textbf{p}=0,  \ x\in \Gamma \}$ with minimal number of boundary conditions such that the boundary term  is never negative, $ \left.\left(\textbf{q}^\top \widetilde{M} \textbf{q}\right)\right|_0^1 \ge 0$.
We will now choose the weight matrix $W$ such that the weighted $L_2$-norm is related to the mechanical energy in the medium. 
Note in particular, if 
\begin{align}
        \label{eq:weight_matrix}
     {W}=\begin{bmatrix}g&0\\0&H\end{bmatrix},
\end{align}
    then the weighted $L_2$-norm is related to the mechanical energy $E$, that is
\begin{align}
        \label{eq:energy}
        \frac{1}{2}\Vert \textbf{q} \Vert_{W}^2=E := \int_{\Omega} \frac12 (gh^2 + Hu^2) \,\text{d}x > 0, \quad \forall \textbf{q} \in \mathbb{R}^2\backslash \{\textbf{0}\}.
    \end{align}
We introduce the boundary term
\begin{align}
        \label{eq:boundary_term}
BT:=-\frac{1}{2gH}\left.\left(\textbf{q}^\top \widetilde{M} \textbf{q}\right)\right|_{0}^{1}= \frac{U}{H} \left( \frac12 h^2 \big\vert_1^0 \right)
        +  \left( uh\big{\vert}_1^0 \right)
        + \frac{U}{g} \left( \frac12 u^2 \big\vert_1^0 \right).
\end{align}
By using the eigen-decomposition of the symmetric matrix $\widetilde{M}$ the boundary term can be re-written as 
  \begin{align}
        \label{eq:boundary_term_0}
BT=\frac{1}{2}\left.\left(\lambda_1 w_1^2 + \lambda_2 w_2^2\right)\right|_{x=0} - \frac{1}{2}\left.\left(\lambda_1 w_1^2 + \lambda_2 w_2^2\right)\right|_{x =1},
\end{align}  
where
\begin{align}\label{eq:w_and_S}
\begin{bmatrix}
            w_1, &w_2
        \end{bmatrix}^\top &= {S}^\top \textbf{q}, \qquad
        {S} =  \begin{bmatrix}
            \frac 1c \left( \lambda_1 - \frac{U}{g}\right) & \frac 1d \left( \lambda_2 - \frac{U}{g}\right) \\
            \frac 1c & \frac{1}{d}
        \end{bmatrix}, 
        \end{align}
        and
    $
        %
        c           = \sqrt{\left( \lambda_1 - \frac{U}{g}\right)^2 + 1},\quad
        d           = \sqrt{\left( \lambda_2 - \frac{U}{g}\right)^2 + 1}.
        $
        Here, ${S}$ is a matrix of orthonormal eigenvectors and  so $S^{\top}S = I$. The eigenvalues, $\lambda_1, \lambda_2$, are real and given by
        {\small
        \begin{align}
        \lambda_1   &= \frac1{2gH} \left(U(g+H) + \sqrt{U^2(g+H)^2 + 4gH(gH - U^2)}\right),\\
        \lambda_2   &= \frac1{2gH} \left(U(g+H) - \sqrt{U^2(g+H)^2 + 4gH(gH - U^2)}\right).
    \end{align}
    }

    The number of boundary conditions will depend on the signs of the eigenvalues $\lambda_1$, $\lambda_2$ which in turn depend on the magnitude of the flow $U$ and the sign of $gH-U^2$. The term $(gH - U^2)$ plays an important role to the change of sign of the eigenvalues. That is, $(gH - U^2)>0$ implies  $\lambda_1 >0$ and $\lambda_2<0$;  $gH-U^2 <0$ implies both of the eigenvalues take the sign of $U$; and $(gH - U^2)=0$ implies one of the eigenvector equals zero, that is $\lambda_1 > 0$, $\lambda_2 =0$ if $U> 0$ and $\lambda_1 = 0$, $\lambda_2 <0$ if $U< 0$. 
    We can also discriminate positive $U>0$ and negative $U<0$. When $U>0$, $x =0$ is an inflow boundary and $x = 1$ is an outflow boundary. The situation reverses when $U<0$, that is, $x =0$ becomes an outflow boundary and $x = 1$ is an inflow boundary.

   \textbf{Sub-critical flow.}
    The flow is sub-critical when $U^2 < gH$ which implies $\lambda_1 >0$ and $\lambda_2 <0$. We need one boundary condition at  $x = 0$ and one boundary condition at $x = 1$. Therefore, for sub-critical flow regime we always need an inflow boundary condition and an outflow boundary condition for any $U$. We formulate the boundary conditions 
    \begin{align}\label{eq:BC-sub-critical}
    \{\mathcal{B}\textbf{p}=\textbf{b},  \ x\in \Gamma \} \equiv \{ {w_1 - \gamma_0 w_2 = b_1(t), \ x =0}; \ {w_2 - \gamma_1 w_1 = b_2(t), \ x = 1}\},
    \end{align}
where $\gamma_0, \gamma_1 \in \mathbb{R}$ are boundary reflection coefficients. The following Lemma constraints the parameters $\gamma_0, \gamma_1$.
\begin{lemma}\label{Lem:BC-Sub-critical flow}
    Consider the boundary term $BT$ defined in \eqref{eq:boundary_term_0} and the boundary condition \eqref{eq:BC-sub-critical} with $\textbf{b} =0$ for sub-critical flows $U^2 < gH$ with $\lambda_1 >0$ and $\lambda_2 <0$. If  $0\leq \gamma_0^2 \leq -\lambda_2/\lambda_1$ and $0\leq \gamma_1^2 \le -\lambda_1/\lambda_2$, then the boundary term is never positive, $BT \le 0$.
\end{lemma}
\begin{proof}
   Let $w_1 = \gamma_0 w_2$ at $ x =0$ and $w_2 = \gamma_1 w_1$ at $ x =1$, and consider
    \begin{align*}
        \left.(\lambda_1 w_1^2 + \lambda_2 w_2^2)\right\vert_0
                - \left.(\lambda_1 w_1^2 + \lambda_2 w_2^2)\right\vert_1                     
            = \left. w_2^2 \left( \lambda_1\gamma_0^2   + \lambda_2  \right) \right\vert_0        
                - \left. w_1^2 \left( \lambda_1  + \lambda_2 \gamma_1^2\right) \right\vert_1.
    \end{align*}
    Thus if $0\leq \gamma_0^2 \leq -\lambda_2/\lambda_1$ and $0\leq \gamma_1^2 \le -\lambda_1/\lambda_2$ then $\left( \lambda_1\gamma_0^2   + \lambda_2  \right) \le 0$ and $\left( \lambda_1  + \lambda_2 \gamma_1^2\right)\ge 0$, and we have
   \begin{align*}
     BT = \frac{1}{2}\left(\left. w_2^2 \left( \lambda_1\gamma_0^2   + \lambda_2  \right) \right\vert_0        
                - \left. w_1^2 \left( \lambda_1  + \lambda_2 \gamma_1^2\right) \right\vert_1\right)\le 0.
  \end{align*}
\end{proof}
    \textbf{Super-critical flow.}
    When  $U^2 > gH$ the flow is super-critical, then $\lambda_1$ and $\lambda_2$ both take the sign of the average flow velocity $U$. That is if $U>0$ then $\lambda_1 > 0$, $\lambda_2 > 0$ and if $U<0$ then $\lambda_1 < 0$, $\lambda_2 < 0$. Thus when $U>0$ we need two boundary conditions at $x =0$ and no boundary conditions at $x =1$. Similarly, when $U<0$ we need two boundary conditions at $x =1$ and no boundary conditions at $x =0$. Therefore, for super-critical flows there are no outflow boundary conditions for any $U$. We formulate the boundary conditions
%
%
  \begin{subequations}\label{eq:BC-super-critical}
        \begin{align}
            \label{eq:BC-super-critical_positive}
    \{\mathcal{B}\textbf{q}=\textbf{b},  \ x\in \Gamma \} \equiv \{ {w_1 = b_1(t), \ w_2 =b_2(t), \ x =0}; \ \text{if $U >0$}\},\\
            \label{eq:BC-super-critical_negative}
    \{\mathcal{B}\textbf{q}=\textbf{b},  \ x\in \Gamma \} \equiv \{ {w_1 = b_1(t), \ w_2 =b_2(t), \ x =1}; \ \text{if $U <0$}\}.
        \end{align}
    \end{subequations}
\begin{lemma}\label{Lem:BC-Super-critical flow}
    Consider the boundary term $BT$ defined in \eqref{eq:boundary_term_0} and the boundary condition \eqref{eq:BC-super-critical} with $\textbf{b} =0$ for super-critical flows $U^2 > gH$, we have $BT \le 0$.
\end{lemma}
\begin{proof}
   Let $U > 0$ with $\lambda_1 > 0$, $\lambda_2 > 0$ if $w_1 = 0, \ w_2 =0$, at  $x =0$, then
    \begin{align*}
        BT = \frac{1}{2}\left(\left.(\lambda_1 w_1^2 + \lambda_2 w_2^2)\right\vert_0
                - \left.(\lambda_1 w_1^2 + \lambda_2 w_2^2)\right\vert_1\right)                     
            = - \frac{1}{2}\left.(\lambda_1 w_1^2 + \lambda_2 w_2^2)\right\vert_1 \le 0.
    \end{align*}
    If $U < 0$ with $\lambda_1 < 0$, $\lambda_2 < 0$  and $w_1 = 0, \ w_2 =0$, at  $x =1$, then we   have
    \begin{align*}
         BT = \frac{1}{2}\left(\left.(\lambda_1 w_1^2 + \lambda_2 w_2^2)\right\vert_0
                - \left.(\lambda_1 w_1^2 + \lambda_2 w_2^2)\right\vert_1\right)                     
            =  \frac{1}{2} \left.(\lambda_1 w^2 + \lambda_2 w_2^2)\right\vert_0 \le 0.
    \end{align*}
\end{proof}  
\textbf{Critical flow.}
The flow is critical when $U^2 = gH$. Note that this case is degenerate, since there is only one  nonzero eigenvalue, that is $U >0$ implies $\lambda_1 > 0$, $\lambda_2 =0$ and $U <0$ implies $\lambda_1 = 0$, $\lambda_2 < 0$.  However, it can also be treated by prescribing only one boundary condition for the system. The location of the boundary condition will be determined by the sign of $U$, similar to the super-critical flow regime. We prescribe the boundary conditions
\begin{subequations}\label{eq:BC-critical}
        \begin{align}
            \label{eq:BC-critical_positive}
    \{\mathcal{B}\textbf{q}=\textbf{b},  \ x\in \Gamma \} \equiv \{ {w_1 = b_1(t), \ x =0}; \ \text{if $U >0$ and $\lambda_2 =0$}\},\\
            \label{eq:BC-critical_negative}
    \{\mathcal{B}\textbf{q}=0,  \ x\in \Gamma \} \equiv \{ w_2 =b_2(t), \ x =1; \ \text{if $U <0$ and $\lambda_1 =0$}\}.
        \end{align}
    \end{subequations}
\begin{lemma}\label{Lem:BC-critical flow}
    Consider the boundary term $BT$ defined in \eqref{eq:boundary_term_0} and the boundary condition \eqref{eq:BC-critical} with $\textbf{b}=0$ for critical flows $U^2 = gH$, we have $BT \le 0$.
\end{lemma}
\begin{proof}
   Let $U > 0$ with $\lambda_1 > 0$, $\lambda_2 = 0$ if $w_1 = 0$, at  $x =0$,
    \begin{align*}
        BT = \frac{1}{2}\left(\left.(\lambda_1 w_1^2 + \lambda_2 w_2^2)\right\vert_0
                - \left.(\lambda_1 w_1^2 + \lambda_2 w_2^2)\right\vert_1\right)                     
            = - \frac{1}{2}\left.\lambda_1 w_1^2 \right\vert_1 \le 0.
    \end{align*}
    If $U < 0$ with $\lambda_1 = 0$, $\lambda_2 < 0$  and $w_2 =0$, at  $x =1$ we also  have
    \begin{align*}
         BT = \frac{1}{2}\left(\left.(\lambda_1 w_1^2 + \lambda_2 w_2^2)\right\vert_0
                - \left.(\lambda_1 w_1^2 + \lambda_2 w_2^2)\right\vert_1\right)                     
            =  \frac{1}{2} \left.\lambda_2 w_2^2\right\vert_0 \le 0.
    \end{align*}
\end{proof}  

We will conclude this section with the theorem which proves the well-posedness of the IBVP \eqref{eq:ibvp}.
\begin{theorem}\label{theo:well-posdeness-ibvp}
    Consider the IBVP \eqref{eq:ibvp} where the boundary operator $\mathcal{B}\textbf{q}=0$ is define by \eqref{eq:BC-sub-critical} with $\gamma_0^2 \leq -\lambda_2/\lambda_1$ and $\gamma_1^2 \le -\lambda_1/\lambda_2$ for sub-critical flows, $U^2 < gH$, by \eqref{eq:BC-super-critical} for the super-critical flow regime, $U^2 > gH$, and by \eqref{eq:BC-critical} for critical flows, $U^2 = gH$, we have the energy estimate
    \begin{align}\label{eq:continuous_energy_estimate}
        \frac{1}{2}\frac{d}{dt}\Vert \textbf{q} \Vert_{W}^2 = gH\times\mathrm{BT}\le 0.
    \end{align}
\end{theorem}
\begin{proof}
    We use the energy method, that is, from the left we multiply \eqref{eq:ibvp_pde} with $\textbf{q}^\top W$ and integrate over the domain. As above integration-by-parts gives 
    $$
\frac{1}{2}\frac{d}{dt}\Vert \textbf{q} \Vert_{W}^2  = \left(\textbf{q},\frac{\partial \textbf{q}}{\partial t}\right)_{W}  = \left(\textbf{q},\textbf{D}\textbf{q}\right)_{W} = gH\times\mathrm{BT}.
    $$
    Using Lemmas \ref{Lem:BC-Sub-critical flow}--\ref{Lem:BC-critical flow} for each flow regime gives $\mathrm{BT}\le 0$, which completes the proof.
    \end{proof}
This energy estimate \eqref{eq:continuous_energy_estimate} is what a stable  numerical method should emulate.
    \begin{table}[b]
        \centering
        \begin{tabular}{|c|c|c|}
        \hline
            Regime         & Type of Boundary &   Number of Boundary condition \\ \hline
            sub-critical   & inflow    & 1          \\
                           & outflow    & 1                             \\ \hline
            critical       & inflow    & 1                           \\
                           & outflow    & 0                          \\ \hline
            super-critical & inflow    & 2                      \\
                           & outflow    & 0                           \\ \hline  
        \end{tabular}
        \caption{The number and location of boundary condition in all regime. The boundary at $x=0$ ($x=1$) is inflow (outflow) boundary if $U>0$ and outflow (inflow) boundary if $U<0$.}
        \label{tab:number-of-boundary-condition}
    \end{table}

\section{Numerical scheme}
\label{sec:num}
    We will now derive a stable  finite volume method for the IBVP \eqref{eq:ibvp} encapsulated in the SBP framework. We will prove numerical stability by deriving discrete energy estimates analogous to Theorem \ref{theo:well-posdeness-ibvp}.
    \subsection{The finite volume method}
    \label{subsec:fvm}
    \begin{figure}[htb!]
        \centering \includegraphics[scale=0.7]{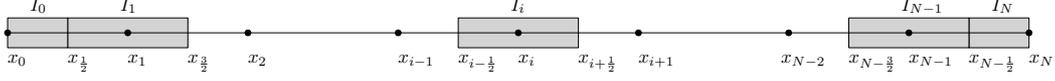}
        \caption{Finite volume nodes $x_i$ and control cells $I_i$.}
        \label{fig:fvm_cells}
    \end{figure}
    To begin, the domain, $\Omega = [0, 1]$, is subdivided into $N+1$ computational nodes having
    $ x_{i} = x_{i-1} + \Delta{x}_i$, for $i = 1, 2, \cdots N$, with $x_0 = 0$, $\Delta{x}_i > 0$ and $\sum_{i=1}^{N}\Delta{x}_i = 1$. We consider the control cell $I_i = [x_{i-\frac12},x_{i+\frac12}]$ for each interior node $1\le i \le N-1$, and for the boundary nodes $\{x_0, x_N\}$ the control cells are $I_{0} = [x_0, x_{1/2}]$ and $I_{N} = [x_{N-1/2}, x_{N}]$, see Figure \ref{fig:fvm_cells}. Note that $|I_i| =\Delta{x}_{i}/2 + \Delta{x}_{i+1}/2$ for the interior nodes $1\le i\le N-1$, and for the boundary nodes $i \in \{0, N\}$ we have $|I_0| = \Delta{x}_{1}/2$ and $|I_N| = \Delta{x}_{N}/2$. The control cells $I_i$ are connected and do not overlap, and $\sum_{i=0}^N |I_i| = \sum_{i=1}^{N}\Delta{x}_i = 1.$

   Consider the integral form of \eqref{eq:ibvp_pde} over the control cells  $I_i$
   {\small
  \begin{subequations}\label{eq:integral_form}
        \begin{align}
            &\dfrac{d}{\text{d}t}\int_{I_0} \textbf{q}(x,t) \,\text{d}x
            + {M} \textbf{q}(x_{\frac12},t) 
                - {M}\textbf{q}(x_{0},t) = 0, \\
       & \dfrac{d}{\text{d}t}\int_{I_i} \textbf{q}(x,t) \,\text{d}x
            + {M} \textbf{q}(x_{i+\frac12},t) 
                - {M}\textbf{q}(x_{i-\frac12},t) = 0,  \quad 1\le i\le N-1,\\                
             &   \dfrac{d}{\text{d}t}\int_{I_N} \textbf{q}(x,t) \,\text{d}x
            + {M} \textbf{q}(x_{N},t) 
                - {M}\textbf{q}(x_{N-1/2},t) = 0 .         
                \label{eq:intgrtd1}
        \end{align}   
         \end{subequations}
         }
     Introduce the cell-average
    \begin{align}
        \bar{\textbf{q}}_i = \frac{1}{|I_i|}\int_{I_i} \textbf{q}(x,t) \text{d}x, 
    \end{align} 
    and approximate the PDE flux $M\textbf{q}$ 
    with the local Lax-Friedrich  flux

    \begin{align}\label{eq:LaxFriedric}
        M\textbf{q}(x_{i+ \frac12},t) \approx \frac{{M}\bar{\textbf{q}}_{i+1} 
            + {M}\bar{\textbf{q}}_i}{2} 
            - \frac{\alpha}2\left(\bar{\textbf{q}}_{i+1}  - \bar{\textbf{q}}_i\right), \quad \alpha \ge 0,
    \end{align} 
    and 
   $
        {M}\textbf{q}(x_{0},t) \approx M \bar{\textbf{q}}_0, \quad {M}\textbf{q}(x_{N},t) \approx M \bar{\textbf{q}}_N.
    $
    The evolution of the cell-average is governed by the semi-dicrete system 
    {\small
     \begin{subequations}\label{eq:semi-discrete}
        \begin{align}
        &|I_0| \frac{d\bar{\textbf{q}_0}}{dt}
            + {M}\frac{\bar{\textbf{q}}_{1} - \bar{\textbf{q}}_{0}}{2} 
            -  \frac{\alpha}{2} (\bar{\textbf{q}}_{1} - \bar{\textbf{q}}_{0}) = 0,\\
         &|I_i| \frac{d\bar{\textbf{q}_i}}{dt} 
            + {M}\frac{\bar{\textbf{q}}_{i+1} - \bar{\textbf{q}}_{i-1}}{2} 
            -  \frac{\alpha}{2} (\bar{\textbf{q}}_{i+1} - 2\bar{\textbf{q}}_{i} + \bar{\textbf{q}}_{i-1}) = 0,  \ 1\le i\le N-1,\\     
          &  |I_N| \frac{d\bar{\textbf{q}_N}}{dt}
            + {M}\frac{\bar{\textbf{q}}_{N}  -\bar{\textbf{q}}_{N-1}}{2} 
            -  \frac{\alpha}{2} (\bar{\textbf{q}}_{N-1} - \bar{\textbf{q}}_{N}) = 0.
        \end{align} 
        \end{subequations}
        }
    Introducing the discrete solution vector $\bar{\textbf{q}} = [\bar{\textbf{q}}_0, \bar{\textbf{q}}_1, \cdots, \bar{\textbf{q}}_N]^{\top}$ and rewriting  \eqref{eq:semi-discrete} in a more compact form, we have  
        \begin{align} \label{eq:num_scheme_1}
        \left(I\otimes{P}\right)\frac{d\bar{\textbf{q}}}{dt} + \left({M}\otimes {Q}\right)  \bar{\textbf{q}} 
                    - \frac{\alpha}{2} \left(I\otimes{A}\right)\bar{\textbf{q}}=0,
        \end{align} 
    where $\otimes$ denotes the Kronecker product and
        \begin{align*}
            {Q} = \begin{pmatrix}
                - \frac12         & \frac12            & 0                            & \cdots & 0        & 0         & 0 \\
                -\frac{1}{2} & 0             & \frac{1}{2}                  & \cdots & 0        & 0         & 0 \\
                \vdots       & \vdots        & \vdots                  & \ddots & \vdots   & \vdots    & \vdots \\
                0                         & 0              & 0               & \cdots & -\frac{1}{2} & 0      & \frac{1}{2} \\
                0                         & 0              & 0               & \cdots & 0        & -\frac12      & \frac12
            \end{pmatrix}, \
            {A} = \begin{pmatrix}
                {-1}           &{1}            & 0                             & \cdots & 0        & 0         & 0 \\
                1           & -2             & 1                   & \cdots & 0        & 0         & 0 \\
                \vdots       & \vdots        & \vdots                  & \ddots & \vdots   & \vdots    & \vdots \\
                0            & 0             & 0                             & \cdots & 1  & -2      & 1 \\
                0            & 0             & 0                             & \cdots & 0        & 1      & -1 
            \end{pmatrix},
        \end{align*}
    and ${P} =  \text{diag} \left([|I_0|,|I_1|,\cdots,|I_N|]\right)$.         
The matrix $Q$ is related to the spatial derivative operator, $A$ is a numerical dissipation operator, and $\alpha \ge 0$ controls the amount of numerical dissipation applied. Note that $A$ is symmetric and negative semi-definite, that is $A = A^\top$ and $ \mathbf{q}^\top A \mathbf{q} \leq 0 $ for all $\mathbf{q} \in \mathbb{R}^{N+1}$.
The important stability property of the semi-discrete approximation \eqref{eq:num_scheme_1} is that the associated discrete derivative operator satisfies the SBP property. To see this, we rewrite 
 equation \eqref{eq:num_scheme_1}  as 
        \begin{align}
            \frac{d\bar{\textbf{q}}}{dt} + \left({M}\otimes {D}_x\right)  \bar{\textbf{q}} 
            - \frac{\alpha}{2} \left(I\otimes{P}^{-1}{A}\right)\bar{\textbf{q}}=0,          
        \end{align}
    where $I$ is the $2\times 2$ identity matrix and
    \begin{align}\label{eq:sbp}
    {D}_x = {P}^{-1}{Q}, \quad 
                {Q}+{Q}^\top = \text{diag}\left([-1,0,\dots,0,1]\right).
        \end{align}
The relation \eqref{eq:sbp} is the so-called SBP property 
\cite{kreiss1974finite,gustafsson2007high} for the first derivative $d/dx$, which can be useful in proving numerical stability of the discrete approximation \eqref{eq:num_scheme_1}. Note that we have not enforced any boundary condition yet, the boundary condition \eqref{eq:ibvp_boundary} will be implemented weakly using penalties.
\subsection{Numerical boundary treatment and  stability}
We will now implement the boundary conditions and prove numerical stability. 
The boundary conditions are implemented using the SAT method, similar terms used as in \cite{carpenter1994time}, by appending the boundary operators \eqref{eq:BC-sub-critical}--\eqref{eq:BC-critical} to the right hand-side of \eqref{eq:num_scheme_1} with penalty weights, we have
        \begin{align} 
            \label{eq:num_scheme_1_SAT}
        \left(I\otimes\textbf{P}\right)\frac{d\bar{\textbf{q}}}{dt} + \left({M}\otimes \textbf{Q}\right)  \bar{\textbf{q}} 
                    - \frac{\alpha}{2} \left(I\otimes\textbf{A}\right)\bar{\textbf{q}}=\mathrm{SAT}.
        \end{align} 
With $\mathbf{e}_0 = [1, 0, \cdots 0]^T$ and $\mathbf{e}_N = [0, 0, \cdots 1]^T$,
the SAT for sub-critical flow is
{\small
\begin{align}\label{eq:sub-critical_SAT}
\mathrm{SAT}=  -\frac{1}{2}\left(W^{-1}SW\otimes\textbf{I}\right)
\begin{bmatrix}
{\tau_{01}} H\mathbf{e}_0\left(\bar{w}_1 - \gamma_0 \bar{w}_2 -b_1(t)\right)
        + { \tau_{N1}} H\mathbf{e}_N \left(\bar{w}_2 - \gamma_1 \bar{w}_1 - b_2(t) \right)\\
{\tau_{02}}g \mathbf{e}_0 \left(\bar{w}_1 - \gamma_0 \bar{w}_2 - b_1(t) \right)
        + {\tau_{N2}} g\mathbf{e}_N\left(\bar{w}_2 - \gamma_1 \bar{w}_1 - b_2(t) \right)
\end{bmatrix},
 \end{align} 
 }
and for critical/super-critical flow regimes we have
{\small
\begin{subequations}\label{eq:super-critical_SAT}
\begin{align} 
\label{eq:positive_U_SAT}
\mathrm{SAT}=  -\frac{1}{2}\left(W^{-1}SW\otimes\textbf{I}\right)
\begin{bmatrix}
    {\tau_{01}} H\mathbf{e}_0(\bar{w}_1-b_1(t)) \\    
         {\tau_{02}}g \mathbf{e}_0 (\bar{w}_2-b_2(t))
\end{bmatrix}, \qquad U > 0,\\
\label{eq:negative_U_SAT}
\mathrm{SAT}=  -\frac{1}{2}\left(W^{-1}SW\otimes\textbf{I}\right)
\begin{bmatrix}
    {\tau_{N1}} H\mathbf{e}_N(\bar{w}_1-b_1(t)) \\    
         {\tau_{N2}}g \mathbf{e}_N (\bar{w}_2-b_2(t))
\end{bmatrix}, \qquad U < 0.
 \end{align} 
 \end{subequations}
 }

\noindent Here $S$ is the eigenvector matrix given in \eqref{eq:w_and_S} and $W$ is the weight matrix defined in \eqref{eq:weight_matrix}. The  coefficients $\tau_{01}$, $\tau_{02}$, $\tau_{N1}$, $\tau_{N2}$ are real penalty parameters to be determined by requiring stability. Note that \eqref{eq:num_scheme_1_SAT}  is a consistent semi-discrete approximation of the IBVP \eqref{eq:ibvp} for all nontrivial choices of the penalty parameters. The semi-discrete approximation \eqref{eq:num_scheme_1_SAT}, given that the discrete derivative operator satisfies the SBP property \eqref{eq:sbp}, is often referred to as the SBP-SAT scheme.
We introduce the discrete weighted $L_2$-norm
$$
\Vert \bar{\textbf{q}} \Vert_{WP}^2:= \bar{\textbf{q}}^T\left(W\otimes P\right)\bar{\textbf{q}}  \ge 0.
$$
Tablehe semi-discrete approximation \eqref{eq:num_scheme_1_SAT} is stable if the discrete energy norm $\Vert \bar{\textbf{q}} \Vert_{WP}^2$ is non-increasing with time.
    We will now prove the stability of the semi-discrete approximation \eqref{eq:num_scheme_1_SAT} for sub-critical flows.
    \begin{theorem}\label{theo:stability sub-critical flow}
        Consider the semi-discrete finite volume approximation \eqref{eq:num_scheme_1_SAT} with the SAT \eqref{eq:sub-critical_SAT} and $\textbf{b} =0$ for sub-critical flow regimes  where $\lambda_1 >0$, $\lambda_2 < 0$ and $\gamma_0^2 \le -\lambda_2/\lambda_1$, $\gamma_1^2 \le -\lambda_1/\lambda_2$. If the penalty parameters are chosen such that 
        $
\tau_{01} = \lambda_1,\quad  \tau_{02} = \gamma_0\lambda_1; \quad \tau_{N2} = -\lambda_2, \quad \tau_{N1} = -\gamma_1\lambda_2,
        $
        then 
        \begin{align*}
        \frac{1}{2}\frac{d}{dt}\Vert \bar{\textbf{q}} \Vert_{WP}^2\le 0, \quad \forall \ t\ge 0.
    \end{align*}
    \end{theorem}
    \begin{proof}
        We use the energy method, that is from the left, we multiply \eqref{eq:num_scheme_1_SAT} with $\bar{\textbf{q}}^T\left(W\otimes \mathbf{I}\right)$ and add the transpose of the product, we have
        {\small
        \begin{align*} 
        \frac{1}{2}\frac{d}{dt}\Vert \bar{\textbf{q}} \Vert_{WP}^2 + \frac{1}{2} \bar{\textbf{q}}^T\left(\widetilde{M}\otimes \left({Q}+{Q}^T\right)\right)  \bar{\textbf{q}} 
                    - \frac{\alpha}{2} \bar{\textbf{q}}^T\left(W\otimes{A}\right)\bar{\textbf{q}}=\bar{\textbf{q}}^T\left(W\otimes \mathbf{I}\right)\mathrm{SAT}.
        \end{align*} 
        }
        Using the SBP property \eqref{eq:sbp} and the eigen-decomposition of $\widetilde{M}$ we have
         {\small
        \begin{align*} 
        \frac{1}{2}\frac{d}{dt}\Vert \bar{\textbf{q}} \Vert_{WP}^2 
                    - \frac{\alpha}{2} \bar{\textbf{q}}^T\left(W\otimes{A}\right)\bar{\textbf{q}}= \frac{1}{2}gH\times\mathrm{BT}_{num},
        \end{align*} 
        }
        where
       \begin{align*}
         \mathrm{BT}_{num} &= \left.\left(\lambda_1 \bar{w}_1^2 + \lambda_2 \bar{w}_2^2 
         -\left(\tau_{01} \bar{w}_1\left(\bar{w}_1 - \gamma_0 \bar{w}_2 \right) + \tau_{02} \bar{w}_2\left(\bar{w}_1 - \gamma_0 \bar{w}_2 \right)\right)\right)\right|_{i=0} \\
         &- \left.\left(\lambda_1 \bar{w}_1^2 + \lambda_2 \bar{w}_2^2 
         +\left(\tau_{N1} \bar{w}_1\left(\bar{w}_2 - \gamma_1 \bar{w}_1 \right) + \tau_{N2} \bar{w}_2\left(\bar{w}_2 - \gamma_1 \bar{w}_1 \right)\right)\right)\right|_{i=N}. 
         \end{align*}
         Thus, if       $
\tau_{01} = \lambda_1,\quad  \tau_{02} = \gamma_0\lambda_1; \qquad \tau_{N2} = -\lambda_2, \quad \tau_{N1} = -\gamma_1\lambda_2,
        $
        then we have
        \begin{align*}
         \mathrm{BT}_{num} &= \left.\left( \lambda_2 +  \lambda_1\gamma_0^2\right)\bar{w}_2^2\right|_{i=0}  - \left.\left( \lambda_1 +  \lambda_2\gamma_1^2\right)\bar{w}_1^2\right|_{i=N}.
         \end{align*}
         Since $\lambda_1 > 0$, $\lambda_2 < 0$ and       
         $$
(\lambda_2 +  \lambda_1\gamma_0^2)\le 0 \iff \gamma_0^2 \le -\lambda_2/\lambda_1; \qquad (\lambda_1 +  \lambda_2\gamma_1^2) \ge 0 \iff \gamma_1^2 \le -\lambda_1/\lambda_2,
        $$
        then we must have $\mathrm{BT}_{num} \le 0$. 
        {Note that for $\alpha \ge 0$ then  $\frac{\alpha}{2} \bar{\textbf{q}}^T\left(W\otimes{A}\right)\bar{\textbf{q}}\leq 0$, and we have}
         {\small
        \begin{align*} 
        \frac{1}{2}\frac{d}{dt}\Vert \bar{\textbf{q}} \Vert_{WP}^2 
                    = \frac{\alpha}{2} \bar{\textbf{q}}^T\left(W\otimes{A}\right)\bar{\textbf{q}}+ \frac{1}{2}gH\times\mathrm{BT}_{num} \le 0.
        \end{align*} 
        }
    \end{proof}
    The next theorem will prove the stability of the semi-discrete approximation \eqref{eq:num_scheme_1_SAT} for super-critical flows.
     \begin{theorem}\label{theo:stability super-critical flow}
        Consider the semi-discrete finite volume approximation \eqref{eq:num_scheme_1_SAT} with the SAT \eqref{eq:super-critical_SAT} and $\textbf{b}=0$ for super-critical flows. If the penalty parameters are chosen such that 
        $
\tau_{01} \ge \lambda_1,\quad  \tau_{02} \ge \lambda_2; \qquad \tau_{N1} \ge -\lambda_1, \quad \tau_{N2} \ge -\lambda_2,
        $
        then 
        \begin{align*}
        \frac{1}{2}\frac{d}{dt}\Vert \bar{\textbf{q}} \Vert_{WP}^2\le 0, \quad \forall \ t\ge 0.
    \end{align*}
    \end{theorem}
    \begin{proof}
    As above the energy method with the SBP property \eqref{eq:sbp} and the eigen-decomposition of $\widetilde{M}$ yields
         {\small
        \begin{align*}
        \frac{1}{2}\frac{d}{dt}\Vert \bar{\textbf{q}} \Vert_{WP}^2- \frac{\alpha}{2} \bar{\textbf{q}}^T\left(W\otimes{A}\right)\bar{\textbf{q}}= \frac{1}{2}gH\times\mathrm{BT}_{num},
        \end{align*} 
        }
        where
        {\small
       \begin{align*}
         &\mathrm{BT}_{num} &= \left.\left((\lambda_1-\tau_{01}) \bar{w}_1^2 + (\lambda_2-\tau_{02}) \bar{w}_2^2 \right)\right|_{i=0} - \left.\left(\lambda_1 \bar{w}_1^2 + \lambda_2 \bar{w}_2^2\right)\right|_{i=N}, \ U> 0,\\
         &\mathrm{BT}_{num} &=  \left.\left(\lambda_1 \bar{w}_1^2 + \lambda_2 \bar{w}_2^2\right)\right|_{i=0} -\left.\left((\lambda_1+\tau_{N1}) \bar{w}_1^2 + (\lambda_2+\tau_{N2}) \bar{w}_2^2 \right)\right|_{i=N}, \ U< 0.
         \end{align*}
         }%
         Therefore, if      
         $ \tau_{01} \ge \lambda_1,  \tau_{02} \ge \lambda_2;  \tau_{N1} \ge -\lambda_1, \tau_{N2} \ge -\lambda_2,$
         then  we have $\mathrm{BT}_{num} \le 0$. Noting that $\alpha \ge 0$ and {as previous, we have $\frac{\alpha}{2} \bar{\textbf{q}}^T\left(W\otimes{A}\right)\bar{\textbf{q}}\leq 0$ which gives us } 
         {
        \begin{align*} 
        \frac{1}{2}\frac{d}{dt}\Vert \bar{\textbf{q}} \Vert_{WP}^2 
                    = \frac{\alpha}{2} \bar{\textbf{q}}^T\left(W\otimes{A}\right)\bar{\textbf{q}}+ \frac{1}{2}gH\times\mathrm{BT}_{num} \le 0.
        \end{align*} 
        }
    \end{proof}
    Finally, we will prove the stability of the semi-discrete approximation \eqref{eq:num_scheme_1_SAT} for critical flows.
     \begin{theorem}\label{theo:stability critical flow}
        Consider the semi-discrete finite volume approximation \eqref{eq:num_scheme_1_SAT} with the SAT \eqref{eq:super-critical_SAT} and $\textbf{b}=0$ for critical flows. If the penalty parameters are chosen such that 
        $
\tau_{01} \ge \lambda_1,\quad  \tau_{02} =0; \quad \tau_{N1} = 0, \quad \tau_{N2} \ge -\lambda_2,
        $
        then 
        \begin{align*}
        \frac{1}{2}\frac{d}{dt}\Vert \bar{\textbf{q}} \Vert_{WP}^2\le 0, \quad \forall \ t\ge 0.
    \end{align*}
    \end{theorem}
    \begin{proof}

    The zero penalties ensure consistency of the SAT, that is $\tau_{02} =0$  and $\tau_{N1} = 0$ give
    {\small
\begin{subequations}
\begin{align*} 
\mathrm{SAT}=  -\frac{1}{2}\left(W^{-1}SW\otimes\textbf{I}\right)
\begin{bmatrix}
    {\tau_{01}} H\mathbf{e}_0\bar{w}_1 \\    
         0
\end{bmatrix}, \qquad U > 0,\\
%
\mathrm{SAT}=  -\frac{1}{2}\left(W^{-1}SW\otimes\textbf{I}\right)
\begin{bmatrix}
    0 \\    
         {\tau_{N2}}g \mathbf{e}_N \bar{w}_2
\end{bmatrix}, \qquad U < 0.
 \end{align*} 
 \end{subequations}
 }
Again the energy method with the SBP property \eqref{eq:sbp} and the eigen-decomposition of $\widetilde{M}$ yield
         {\small
        \begin{align*}
        \frac{1}{2}\frac{d}{dt}\Vert \bar{\textbf{q}} \Vert_{WP}^2- \frac{\alpha}{2} \bar{\textbf{q}}^T\left(W\otimes{A}\right)\bar{\textbf{q}}= \frac{1}{2}gH\times\mathrm{BT}_{num},
        \end{align*} 
        }
        where
        {\small
       \begin{align*}
         &\mathrm{BT}_{num} = \left.(\lambda_1-\tau_{01}) \bar{w}_1^2\right|_{i=0} - \left.\lambda_1 \bar{w}_1^2 \right|_{i=N}, \ U> 0, \ \lambda_1 >0, \ \lambda_2 =0 \\
         &\mathrm{BT}_{num} =  \left. \lambda_2 \bar{w}_2^2\right|_{i=0} -\left.(\lambda_2+\tau_{N2}) \bar{w}_2^2\right|_{i=N}, \ U< 0,  \ \lambda_1 =0, \ \lambda_2 <0.
         \end{align*}
         }
         Therefore, if      
         $
\tau_{01} \ge \lambda_1,$ and $ \tau_{N2} \ge -\lambda_2,
        $
        then we have $\mathrm{BT}_{num} \le 0$. 
        {Using the fact that $\alpha \ge 0$ and $\frac{\alpha}{2} \bar{\textbf{q}}^T\left(W\otimes{A}\right)\bar{\textbf{q}}\leq 0$ again} gives 
         {\small
        \begin{align*} 
        \frac{1}{2}\frac{d}{dt}\Vert \bar{\textbf{q}} \Vert_{WP}^2 
                    = \frac{\alpha}{2} \bar{\textbf{q}}^T\left(W\otimes{A}\right)\bar{\textbf{q}}+ \frac{1}{2}gH\times\mathrm{BT}_{num} \le 0.
        \end{align*} 
        }
    \end{proof}

\section{Numerical experiments}
    \label{sec:num-exp}

    In this section, we perform numerical experiments  to verify the analysis undertaken in the previous sections.
    Similar to the theoretical analysis, the numerical experiments cover the three flow regimes, namely sub-critical, critical and super-critical flow regimes.
    We used $H=1$ m,  $g=9.8$ m/s$^2$, and  $U \in \{\frac 12 \sqrt{gH},  \sqrt{gH}, 2 \sqrt{gH}\}$, which correspond to the three different flow regimes. The interval of interest is $[0,L]$ with $L> 0$. Note that $U>0$ so that $x=0$ is the in-flow boundary and $x = L$ is the outflow boundary. The locations and the number of boundary conditions required are given in Table \ref{tab:number-of-boundary-condition}, and the explicit forms of the boundary conditions considered here are given in Table \ref{table:boundarycondition-setup} where $g_1(t)$ and $g_2(t)$ are the boundary data. 
    \begin{table}[h!]
        \centering
        \begin{tabular}{|l|c|c|c|c|}
            \hline
            Regime         & $ U $   & Boundaries &  {Boundary conditions}  \\
            \hline
            sub-critical   & $< \sqrt{gH}$      & $x = 0$ & $\frac12(h + \sqrt{H/g}\,u) =g_1$ \\ 

                           &                     & $x = L$ & $\frac12(h - \sqrt{H/g}\, u) =g_2$ \\ \hline

            critical       & $=\sqrt{gH}$         & $x = 0$ & $\frac12(h + \sqrt{H/g}\,u) =g_1$ \\ \hline 
            super-critical   
            & $>\sqrt{gH}$  & $x = 0$ &$\frac12(h +\sqrt{H/g}\, u) =g_1$  \\
            &               &         & $\frac12 (h - \sqrt{H/g}\, u) =g_2$ \\ \hline 
        \end{tabular}
        \caption{Transmissive boundary conditions in all regimes with $U>0$.}
        \label{table:boundarycondition-setup}
    \end{table}

    The semi-discrete system \eqref{eq:num_scheme_1_SAT} is integrated in time using the classical fourth-order accurate explicit Runge-Kutta method with the  time step 
    \begin{align*}
        \Delta t = \mathrm{Cr} \, \frac{\Delta x}{|U|+\sqrt{gH}}, \quad  \mathrm{Cr} = 0.25,
    \end{align*}
    where $\mathrm{Cr}$ is the Courant–Friedrichs–Lewy number, $\Delta x = L/N$ is the uniform cell width and $N$ is the number of finite volume cells. We will consider a centered numerical flux with $\alpha =0$  and the local Lax-Friedrich's numerical flux \eqref{eq:LaxFriedric} with 
    $\alpha>0$,  and verify numerical accuracy.  Note that the semi-discrete approximation is energy stable for all $\alpha  \ge 0$.

    \paragraph{Non-homogeneous boundary data.} 

    We consider zero initial conditions, that is $u(x,0)=0$ and $h(x,0)=0$, and send a wave into the domain through the in-flow boundary at $x=0$. We will consider specifically $g_1(t)\ne 0$ and $g_2(t)=0$, for the boundary conditions given in Table \ref{table:boundarycondition-setup}, so that the corresponding IBVP has the exact solution
$$ \label{eq:exact_solution}
h(x, t) = g_1\left(t-\frac{x}{U+\sqrt{gH}}\right), \quad u(x, t) = \frac{1}{\sqrt{H/g}}g_1\left(t-\frac{x}{U+\sqrt{gH}}\right).
$$
We will consider a smooth boundary data given by

    \begin{equation}\label{eq:boundary_data_smooth_pulse}
        g_1(t)= \begin{cases}(\sin(\pi t))^ 4 & \text { if } 0\leq t \leq 1,  \\ 0 & \text { otherwise},\end{cases}
        \qquad g_2(t) = 0, \quad \forall \ t\ge 0,
    \end{equation}
and a non-smooth boundary data given by
\begin{equation}\label{eq:boundary_data_stepfunct}
    g_1(t)= \begin{cases}1 & \text { if } 0< t \leq 1,  \\ 0 & \text { otherwise},\end{cases}
    \qquad g_2(t) = 0, \quad \forall \ t\ge 0.
    \end{equation}
    The boundary data  for the boundary conditions in Table \ref{table:boundarycondition-setup} can be rewritten as $b_1(t)$ and/or $b_2(t)$ and in terms of $w_1$ and $w_2$ for the given boundary. 
    
    In the sub-critical case, by utilizing the identity \eqref{eq:w_and_S}, the boundary condition can be rewritten in the form \eqref{eq:BC-sub-critical} with
    $$
    \gamma_0 =  - \frac{
        \frac{1}{d}\left( \sqrt[]{\frac gH} \left(\lambda_2 -\frac{ U}{g}\right)-1 \right)
    }
    {
        \frac{1}{c}\left( \sqrt[]{\frac gH} \left(\lambda_1 -\frac{ U}{g}\right)-1 \right)
    }
    , \quad 
    \gamma_1 =  - \frac{
        \frac{1}{d}\left( \sqrt[]{\frac gH} \left(\lambda_2 -\frac{ U}{g}\right)+1 \right)
    }
    {
        \frac{1}{c}\left( \sqrt[]{\frac gH} \left(\lambda_1 -\frac{ U}{g}\right)+1 \right)
    }.
    $$
    We can show that $\gamma_0$ and $\gamma_1$ satisfy the condition of Lemma \ref{Lem:BC-Sub-critical flow}, that is $\gamma_0^2 \leq -\lambda_2/\lambda_1$ and $\gamma_1^2 \le -\lambda_1/\lambda_2$ for all $|U|<\sqrt{gH}$.

    For the critical flow regime, we have $U^2=gH$ and $\lambda_2=0$. Only one boundary condition is imposed at the inflow, $\bar{w_1} = b_1(t)$. This condition can be rewritten to match the condition in Theorem \ref{theo:stability critical flow} by using  the identity \eqref{eq:w_and_S} and the fact that $U^2=gH$. 

    For the super-critical flow regime, two boundary conditions are imposed at the inflow boundary. That is $\bar{w}_1 = b_1(t)$, $\bar{w}_2=b_2(t)\equiv 0$ as the boundary condition at $x = 0$. This conditions is equivalent to \eqref{eq:BC-super-critical_positive}.

The boundary data will generate a pulse from the left boundary at $x =0$, which will propagate through the domain and leave the domain through $x = L$.

\begin{figure}[h!]
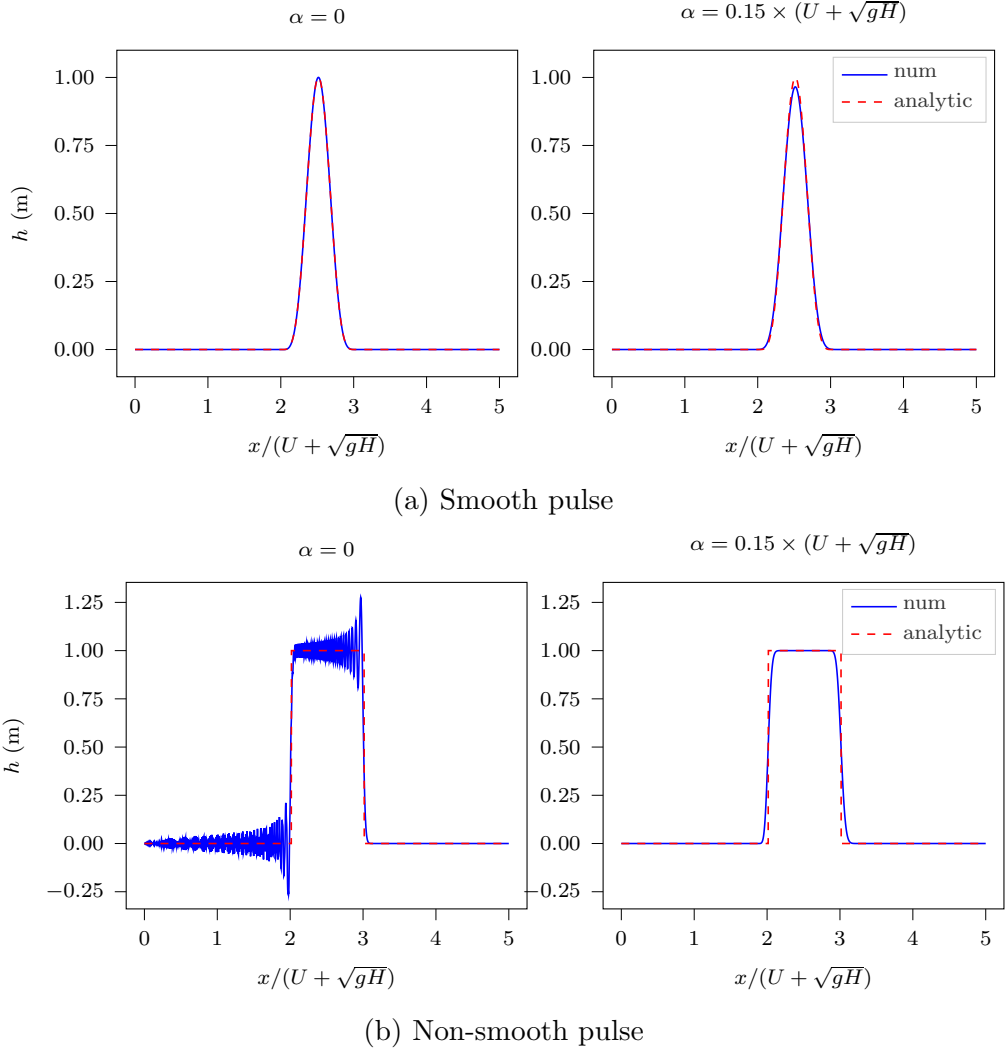

    \centering
    \subcaptionbox{Smooth pulse}{
        \input{numresult/v7.2.1-gaussian-subcritical-dissipation-vary}
        \label{fig:compare-dissipation-gaussian}
    }
    \subcaptionbox{Non-smooth pulse}{
    \input{numresult/v7.2.1-stepfunct-subcritical-dissipation-vary}
    \label{fig:compare-dissipation-stepfunction}
    }
    
    \caption{The snapshots of the numerical and exact solutions with $\Delta x = L\times 2^{-11}$ m at time $t=3.02$ s for a sub-critical flow regime with smooth and non-smooth boundary data.  For the smooth boundary data the numerical solution matches the exact solution well for $\alpha = 0$ and $\alpha = 0.15 \times (U+\sqrt{gH})>0$. Note, however, with $\alpha = 0.15 \times (U+\sqrt{gH})>0$ the peak of the numerical solution is slightly dissipated. For the non-smooth boundary data, when $\alpha=0$, the propagation speed of pulse is well approximated by the numerical solution. However, there are numerical oscillations generated by the propagating discontinuities. When $\alpha = 0.15 \times (U+\sqrt{gH})>0$ the numerical solution is non-oscillatory, but propagating the discontinuities  are smoothed out. 
    }
\label{fig:flow-snapshot}
\end{figure}

 Figure \ref{fig:flow-snapshot} shows the snapshot of the sub-critical flow solutions at $t = 3.02$ s for both smooth and non-smooth boundary data, with $\alpha = 0$ and $\alpha = 0.15 \times (U+\sqrt{gH})>0$. In the plots, we have scaled the horizontal axis by the wave speed $(U+\sqrt{gH})$ so that the solution is spatially invariant for all flow regimes. Note that for the smooth pulse the numerical solution matches the exact solution excellently well for $\alpha = 0$ and $\alpha = 0.15 \times (U+\sqrt{gH})>0$. Although with $\alpha = 0.15 \times (U+\sqrt{gH})>0$ the peak of the numerical slightly dissipated. For the non-smooth pulse, when $\alpha=0$, the propagation speed of pulse is well approximated by the numerical solution. However, there are numerical oscillations generated by the propagating discontinuities. When $\alpha = 0.15 \times (U+\sqrt{gH})>0$ the numerical solution is non-oscillatory, but the discontinuous edges of the solutions are smeared. 

\begin{figure}[h!]
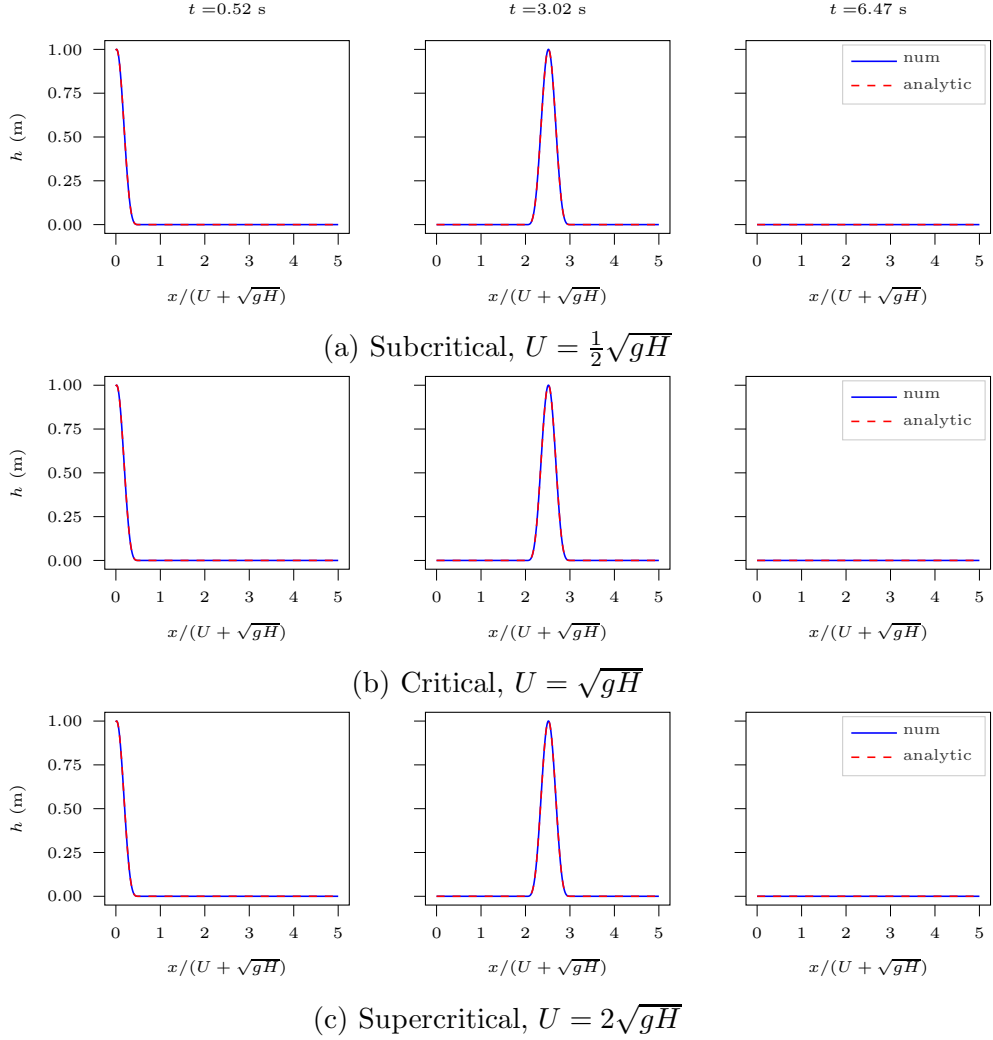

    \centering
    \subcaptionbox{Subcritical, $U=\frac 12 \sqrt{gH}$}{        
        \input{numresult/v7.2.1-gaussian-subcritical-dissipation-0.0-finaltime-6.483185307179586}
    }
    \subcaptionbox{Critical, $U=\sqrt{gH}$}{
        \input{numresult/v7.2.1-gaussian-critical-dissipation-0.0-finaltime-6.483185307179586}
    }
    \subcaptionbox{Supercritical, $U=2 \sqrt{gH}$}{
        \input{numresult/v7.2.1-gaussian-supercritical-dissipation-0.0-finaltime-6.483185307179586}
    }
    \caption{The evolution of the numerical solutions and the exact solutions, at all the three flow regimes with smooth boundary data, $\Delta x = L\times 2^{-11}$ m and $\alpha =0$. The solutions enter the domain through the in-flow boundary at $x=0$ and leave the domain through the out-flow  at $x=L= (U+\sqrt{gH})\times5$.  Note that because of the re-scaling of  the $x$-axis to $x/(U+\sqrt{gH})$, the solutions are invariant for the all three flow regimes. 
            }
\label{fig:flow-evolution-gaussian}
\end{figure}

   The evolution of the numerical solutions and the exact solutions, at all flow regimes  are shown in Figure \ref{fig:flow-evolution-gaussian} for the smooth pulse and in Figure \ref{fig:flow-evolution-step-function}  for the non-smooth pulse. The pulses enter the domain through the in-flow boundary at $x=0$ and leave the domain through the out-flow  at $x=L= (U+\sqrt{gH})\times5$.  Note that because of the re-scaling of  the $x$-axis to $x/(U+\sqrt{gH})$, the solutions are invariant for all three flow regimes. 

\begin{figure}[h!]
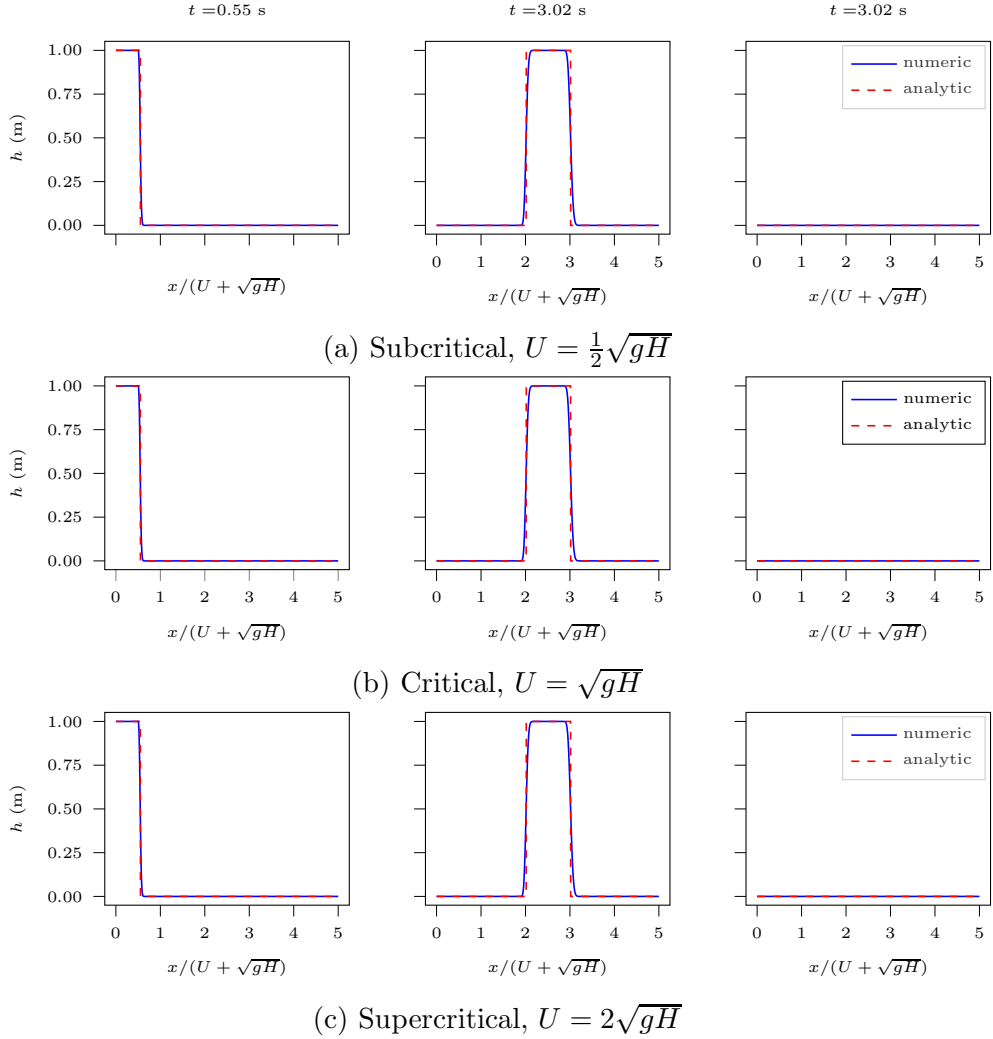

    \centering
    \label{fig:flow-evolution-step-function}
    \subcaptionbox{Subcritical, $U=\frac 12 \sqrt{gH}$}{
        \input{numresult/v7.2.1-stepfunct-subcritical-dissipation-0.7043614129124337-finaltime-6.483185307179586}
    }
    \subcaptionbox{Critical, $U=\sqrt{gH}$}{
        \input{numresult/v7.2.1-stepfunct-critical-dissipation-0.9391485505499116-finaltime-6.483185307179586}
    }
    
    \subcaptionbox{Supercritical, $U=2\sqrt{gH}$}{
        \input{numresult/v7.2.1-stepfunct-supercritical-dissipation-1.4087228258248674-finaltime-6.483185307179586}
    }     
    \caption{ The evolution of the numerical solutions and the exact solutions, at all three flow regimes with non-smooth boundary data, $\Delta x = L\times 2^{-11}$ m and $\alpha =0.15\times(U+\sqrt{gH})>0$. The discontinuous solutions enter the domain through the in-flow boundary at $x=0$ and leave the domain through the out-flow  at $x=L= (U+\sqrt{gH})\times5$.  Note that because of the re-scaling of  the $x$-axis to $x/(U+\sqrt{gH})$, the solutions are invariant for all the three flow regimes. 
            }
 \label{fig:flow-evolution-step-function}
\end{figure}

\begin{figure}[h!]
        \centering
        \subcaptionbox{Subcritical,  ${U} = \frac 12 \sqrt{gH}$}{%
            \input{numresult/v2-transmissive-dissipative-subcritical-error.tex}%
        }
        \subcaptionbox{Critical,  ${U} =\sqrt{gH}$}{%
            \input{numresult/v2-transmissive-dissipative-critical-error.tex}%
        }   
        \subcaptionbox{Supercritical,  ${U} = 2 \sqrt{gH}$}{%
            \input{numresult/v2-transmissive-dissipative-supercritical-error.tex}%
        }
        \caption{The error and  convergence of the error at final time $t=0.1$ using manufactured solution  for all flow regimes.}
            \label{fig:sub-critical-dissipative-convergencetest}
        \label{fig:convergence-test-error}
    \end{figure}
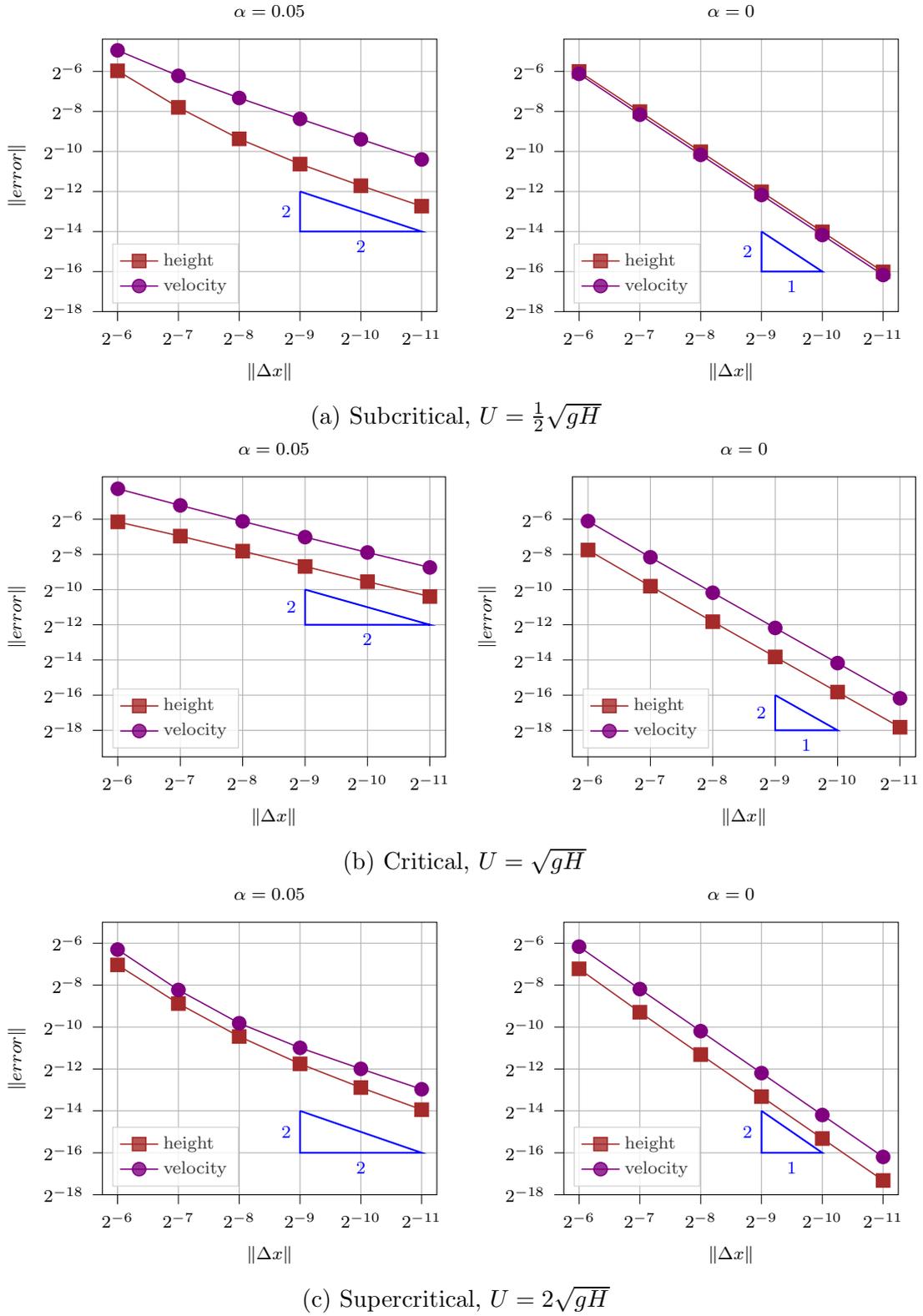

\paragraph{Convergence test.}
    Here, we verify the convergence properties of the numerical method. We will use the method of manufactured solution \cite{roache_code_2001}. That is, we force the system to have the exact smooth solution 
    \begin{align}
        \label{eq:mms}
        h(x,t) = \cos(2\pi t) \sin(6 \pi x ), \qquad 
        u(x,t) = \sin(2\pi t) \cos(4 \pi x ).
    \end{align}
    The initial conditions $h(x,0)$ $u(x,0)$ and the boundary data $g_1(t)$ and $g_2(t)$ are chosen to match the analytical solution \eqref{eq:mms}.
    We compute the numerical solution on a sequence of increasing number of finite volume cells, 
    $N=64,128,256,512,1024,2048$.  
    The $L_2$-error and convergence rates of the error are shown in Figure \ref{fig:convergence-test-error} and also presented in Table \ref{table:convergence-test-error}. We have performed numerical experiments with no dissipation $\alpha = 0$ and with numerical dissipation set on $\alpha = 0.05$. From Table \ref{table:convergence-test-error} we see that the method is second order accurate $O(\Delta{x}^2)$ when $\alpha =0$ and first order accurate $O(\Delta{x})$ when $\alpha >0$. These are in agreement with the theory.
    
    \begin{table}[h!]
    \centering
        \caption{The error and  convergence of the error at final time $t=0.1$ using manufactured solution  for all flow regimes.}
        \label{table:convergence-test-error}
        \subcaptionbox{Subcritical,  ${U} = \frac 12 \sqrt{gH}$}
        {\input{numresult/error_table_subcritical_v3}}
        \\[6pt]
        \subcaptionbox{Critical,  ${U} =\sqrt{gH}$}
        {\input{numresult/error_table_critical_v3}}
        \\[6pt]
        \subcaptionbox{Supercritical,  ${U} = 2 \sqrt{gH}$}
        {\input{numresult/error_table_supercritical_v3}}
    \end{table}
    
\section{Conclusion}
    \label{sec:conc}  
    Well-posed boundary conditions are crucial for accurate numerical solutions of IBVPs. In this study we have analysed well-posed boundary conditions for the linear SWWE in 1D.
    The analysis is based on the energy method and prescribes the number, location and form of the boundary conditions so that the IBVP is well-posed. A summary of the result are shown in the Table \ref{tab:number-of-boundary-condition}, and covers all flow regimes. We formulate the boundary conditions such that they can readily implemented in a stable manner using the SBP-SAT method. We propose a finite volume method formulated in SBP framework and implement the boundary conditions weakly using SAT. 
    Stable penalty parameters and prove of numerical stability derived via discrete energy estimates analogous to the continuous estimate. Numerical experiments are performed to verify the analysis. The error rates comply with the methods that we use. 
    Our continuous and numerical analysis covers all flow regimes, and can be extended to the nonlinear problem. The next step in our study will extend the 1D theory and results to 2D, and implement our scheme in open source software \cite{Roberts_ANUGA_2022, nielsen2005hydrodynamic} for efficient and accurate simulations of the nonlinear shallow water equations.

\paragraph{Acknowledgements}
This research is conducted as part of doctoral study funded by Indonesian Endowment Fund for Education (LPDP).

\ifx\printbibliography\undefined
    \bibliographystyle{plain}
    \bibliography{references}
\else\printbibliography\fi

\end{document}

%% file: numresult/error_table_subcritical_v3.tex
\begin{tabular}{lrrrrrrrr}
\hline
$N$  & \multicolumn{4}{c}{$\alpha=0.0$} &  \multicolumn{4}{c}{$\alpha=0.05$}   \\ \hline
     & $h$ error & rate & $u$ error & rate 
     & $h$ error & rate & $u$ error & rate  \\ \hline
64   &  1.56e-02 &          &  1.44e-02 &          &  1.60e-02 &          &  3.24e-02 &          \\
128  &  3.88e-03 &        2.01 &  3.49e-03 &        2.06 &  4.50e-03 &        1.77 &  1.34e-02 &        1.21 \\
256  &  9.68e-04 &        2.00 &  8.69e-04 &        2.01 &  1.51e-03 &        1.49 &  6.23e-03 &        1.08 \\
512  &  2.42e-04 &        2.00 &  2.17e-04 &        2.01 &  6.31e-04 &        1.20 &  3.02e-03 &        1.03 \\
1024 &  6.05e-05 &        2.00 &  5.41e-05 &        2.00 &  2.98e-04 &        1.06 &  1.49e-03 &        1.01 \\
2048 &  1.51e-05 &        2.00 &  1.35e-05 &        2.00 &  1.47e-04 &        1.01 &  7.41e-04 &        1.01 \\
\hline
\end{tabular}

%% file: numresult/error_table_critical_v3.tex
\begin{tabular}{lrrrrrrrr}
\hline
$N$  & \multicolumn{4}{c}{$\alpha=0.0$} &  \multicolumn{4}{c}{$\alpha=0.05$}   \\ \hline
     & $h$ error & rate & $u$ error & rate 
     & $h$ error & rate & $u$ error & rate  \\ \hline
64   &  4.64e-03 &           &  1.45e-02 &           &  1.41e-02 &           &  5.19e-02 &           \\
128  &  1.12e-03 &        2.08 &  3.50e-03 &        2.08 &  8.03e-03 &        0.88 &  2.70e-02 &        0.96 \\
256  &  2.76e-04 &        2.03 &  8.63e-04 &        2.03 &  4.44e-03 &        0.90 &  1.44e-02 &        0.94 \\
512  &  6.88e-05 &        2.00 &  2.15e-04 &        2.00 &  2.43e-03 &        0.92 &  7.72e-03 &        0.93 \\
1024 &  1.72e-05 &        2.00 &  5.39e-05 &        2.00 &  1.33e-03 &        0.91 &  4.21e-03 &        0.92 \\
2048 &  4.30e-06 &        2.00 &  1.35e-05 &        2.00 &  7.43e-04 &        0.90 &  2.34e-03 &        0.90 \\
\hline
\end{tabular}

%% file: numresult/error_table_supercritical_v3.tex
\begin{tabular}{lrrrrrrrr}
\hline
$N$  & \multicolumn{4}{c}{$\alpha=0.0$} &  \multicolumn{4}{c}{$\alpha=0.05$}   \\ \hline
     & $h$ error & rate & $u$ error & rate 
     & $h$ error & rate & $u$ error & rate  \\ \hline
64   &  6.71e-03 &           &  1.40e-02 &           &  7.63e-03 &           &  1.27e-02 &           \\
128  &  1.60e-03 &        2.10 &  3.43e-03 &        2.03 &  2.12e-03 &        1.80 &  3.33e-03 &        1.90 \\
256  &  3.93e-04 &        2.03 &  8.54e-04 &        2.01 &  7.15e-04 &        1.49 &  1.11e-03 &        1.50 \\
512  &  9.79e-05 &        2.01 &  2.13e-04 &        2.00 &  2.90e-04 &        1.23 &  4.91e-04 &        1.13 \\
1024 &  2.45e-05 &        2.00 &  5.32e-05 &        2.00 &  1.32e-04 &        1.10 &  2.45e-04 &        1.00 \\
2048 &  6.11e-06 &        2.00 &  1.33e-05 &        2.00 &  6.35e-05 &        1.04 &  1.25e-04 &        0.98 \\
\hline
\end{tabular}

%% file: numresult/v2-transmissive-dissipative-subcritical-error.tex
\begin{tikzpicture}[font=\scriptsize]

\definecolor{brown}{RGB}{165,42,42}
\definecolor{darkgray176}{RGB}{176,176,176}
\definecolor{lightgray204}{RGB}{204,204,204}
\definecolor{purple}{RGB}{128,0,128}

\begin{groupplot}[group style={group size=2 by 1,horizontal sep= 2cm},width=0.5*\textwidth]
\nextgroupplot[
legend cell align={left},
legend style={fill opacity=0.8, draw opacity=1, text opacity=1, draw=lightgray204},
legend pos = south west,
tick align=outside,
tick pos=left,
title={$\alpha = 0.05$},
x dir=reverse,
x grid style={darkgray176},
xlabel={\(\displaystyle   \Vert \Delta x \Vert \)},
xmajorgrids,
xmin=-11.25, xmax=-5.75,
xtick style={color=black},
y grid style={darkgray176},
ylabel={\(\displaystyle   \Vert error \Vert \)},
ymajorgrids,
xtick={-6,-7,-8,-9,-10,-11},
xticklabels={$2^{-6}$,$2^{-7}$,$2^{-8}$,$2^{-9}$,$2^{-10}$,$2^{-11}$,},
ytick={-18,-16,-14,-12,-10,-8,-6},
yticklabels={$2^{-18}$,$2^{-16}$,$2^{-14}$,$2^{-12}$,$2^{-10}$,$2^{-8}$,$2^{-6}$},
ymin=-18, ymax=-4.3860557387139,
ytick style={color=black},
y post scale=1
]
\addplot [semithick, brown, mark=square*, mark size=3, mark options={solid}]
table {%
-6 -5.96776805779097
-7 -7.79525462860847
-8 -9.36931541795807
-9 -10.629591324947
-10 -11.7119970319054
-11 -12.7334562099269
};
\addlegendentry{height}
\addplot [semithick, purple, mark=*, mark size=3, mark options={solid}]
table {%
-6 -4.94734028891731
-7 -6.21977207711996
-8 -7.3261489766236
-9 -8.37099564532772
-10 -9.39010694812238
-11 -10.3986318094247
};
\addlegendentry{velocity}
\addplot [thick, blue]
table {%
-9 -12
-9 -14
-11 -14
-9 -12
};
\node[left,blue] at (-9,-13) {2};
\node[below,blue] at (-10,-14) {2};

\nextgroupplot[
legend cell align={left},
legend style={fill opacity=0.8, draw opacity=1, text opacity=1, draw=lightgray204},
legend pos = south west,
scaled y ticks=manual:{}{\pgfmathparse{#1}},
tick align=outside,
tick pos=left,
title={$\alpha = 0$},
x dir=reverse,
x grid style={darkgray176},
xlabel={\(\displaystyle \Vert \Delta x \Vert \)},
xmajorgrids,
xmin=-11.25, xmax=-5.75,
xtick style={color=black},
y grid style={darkgray176},
ymajorgrids,
ytick distance=2,
ymin=-18, ymax=-4.3860557387139,
ytick style={color=black},
xtick={-6,-7,-8,-9,-10,-11},
xticklabels={$2^{-6}$,$2^{-7}$,$2^{-8}$,$2^{-9}$,$2^{-10}$,$2^{-11}$,},
ytick={-18,-16,-14,-12,-10,-8,-6},
yticklabels={$2^{-18}$,$2^{-16}$,$2^{-14}$,$2^{-12}$,$2^{-10}$,$2^{-8}$,$2^{-6}$},
y post scale=1,
]

\addplot [semithick, brown, mark=square*, mark size=3, mark options={solid}]
table {%
-6 -6.00185695717284
-7 -8.00918732269432
-8 -10.0122412035763
-9 -12.0124852110179
-10 -14.0125547624292
-11 -16.0128437302733
};
\addlegendentry{height}
\addplot [semithick, purple, mark=*, mark size=3, mark options={solid}]
table {%
-6 -6.12067174263269
-7 -8.16317037660527
-8 -10.1682015908377
-9 -12.172783073616
-10 -14.1740978290396
-11 -16.1730312929854
};
\addlegendentry{velocity}

\addplot [thick, blue]
table {%
-9 -14
-9 -16
-10 -16
-9 -14
};
\node[left,blue] at (-9,-15) {2};
\node[below,blue] at (-9.5,-16) {1};
\end{groupplot}

\end{tikzpicture}

%% file: numresult/v2-transmissive-dissipative-critical-error.tex
\begin{tikzpicture}[font=\scriptsize]

\definecolor{brown}{RGB}{165,42,42}
\definecolor{darkgray176}{RGB}{176,176,176}
\definecolor{lightgray204}{RGB}{204,204,204}
\definecolor{purple}{RGB}{128,0,128}

\begin{groupplot}[group style={group size=2 by 1, horizontal sep= 2cm}, width=0.51*\textwidth]
\nextgroupplot[
legend cell align={left},
legend style={fill opacity=0.8, draw opacity=1, text opacity=1, draw=lightgray204},
legend pos = south west,
tick align=outside,
tick pos=left,
title={$\alpha = 0.05$},
x dir=reverse,
x grid style={darkgray176},
xlabel={\(\displaystyle \Vert \Delta x \Vert \)},
xmajorgrids,
xmin=-11.25, xmax=-5.75,
xtick style={color=black},
y grid style={darkgray176},
ylabel={\(\displaystyle  \Vert error \Vert \)},
ymajorgrids,
xtick={-6,-7,-8,-9,-10,-11},
xticklabels={$2^{-6}$,$2^{-7}$,$2^{-8}$,$2^{-9}$,$2^{-10}$,$2^{-11}$,},
ytick={-18,-16,-14,-12,-10,-8,-6},
yticklabels={$2^{-18}$,$2^{-16}$,$2^{-14}$,$2^{-12}$,$2^{-10}$,$2^{-8}$,$2^{-6}$},
ymin=-19.5044826877738, ymax=-3.58913217258066,
ytick style={color=black},
y post scale=1
]
\addplot [semithick, brown, mark=square*, mark size=3, mark options={solid}]
table {%
-6 -6.14938775078
-7 -6.95997886995496
-8 -7.81364639075632
-9 -8.68610782933716
-10 -9.54920217254814
-11 -10.3948179130136
};
\addlegendentry{height}
\addplot [semithick, purple, mark=*, mark size=3, mark options={solid}]
table {%
-6 -4.26710265054399
-7 -5.21173794135273
-8 -6.12120521425776
-9 -7.01680557222345
-10 -7.89045750426652
-11 -8.74108385745008
};
\addlegendentry{velocity}
\addplot [thick, blue]
table {%
-9 -10
-9 -12
-11 -12
-9 -10
};
\node[left,blue] at (-9,-11) {2};
\node[below,blue] at (-10,-12) {2};

\nextgroupplot[
legend cell align={left},
legend style={fill opacity=0.8, draw opacity=1, text opacity=1, draw=lightgray204},
legend pos = south west,
tick align=outside,
tick pos=left,
title={$\alpha = 0$},
x dir=reverse,
x grid style={darkgray176},
xlabel={\(\displaystyle \Vert \Delta x \Vert \)},
xmajorgrids,
xmin=-11.25, xmax=-5.75,
xtick style={color=black},
y grid style={darkgray176},
ylabel={\(\displaystyle  \Vert error \Vert \)},
ymajorgrids,
xtick={-6,-7,-8,-9,-10,-11},
xticklabels={$2^{-6}$,$2^{-7}$,$2^{-8}$,$2^{-9}$,$2^{-10}$,$2^{-11}$,},
ytick={-18,-16,-14,-12,-10,-8,-6},
yticklabels={$2^{-18}$,$2^{-16}$,$2^{-14}$,$2^{-12}$,$2^{-10}$,$2^{-8}$,$2^{-6}$},
ymin=-19.5044826877738, ymax=-3.58913217258066,
ytick style={color=black},
y post scale=1
]
\addplot [semithick, brown, mark=square*, mark size=3, mark options={solid}]
table {%
-6 -7.75029629082354
-7 -9.80679954434251
-8 -11.825532227065
-9 -13.8264692389388
-10 -15.8259481424636
-11 -17.8265122098105
};
\addlegendentry{height}
\addplot [semithick, purple, mark=*, mark size=3, mark options={solid}]
table {%
-6 -6.10390541620671
-7 -8.1604086697293
-8 -10.179141352442
-9 -12.1800783642531
-10 -14.1795572675144
-11 -16.1801213338056
};
\addlegendentry{velocity}
\addplot [thick, blue]
table {%
-9 -16
-9 -18
-10 -18
-9 -16
};
\node[left,blue] at (-9,-17) {2};
\node[below,blue] at (-9.5,-18) {1};
\end{groupplot}

\end{tikzpicture}

%% file: numresult/v2-transmissive-dissipative-supercritical-error.tex
\begin{tikzpicture}[font=\scriptsize]

\definecolor{brown}{RGB}{165,42,42}
\definecolor{darkgray176}{RGB}{176,176,176}
\definecolor{lightgray204}{RGB}{204,204,204}
\definecolor{purple}{RGB}{128,0,128}

\begin{groupplot}[group style={group size=2 by 1,horizontal sep= 2cm},width=0.5*\textwidth]
\nextgroupplot[
legend cell align={left},
legend style={fill opacity=0.8, draw opacity=1, text opacity=1, draw=lightgray204},
legend pos = south west,
tick align=outside,
tick pos=left,
title={$\alpha = 0.05$},
x dir=reverse,
x grid style={darkgray176},
xlabel={\(\displaystyle  \Vert \Delta x \Vert \)},
xmajorgrids,
xmin=-11.25, xmax=-5.75,
xtick style={color=black},
y grid style={darkgray176},
ylabel={\(\displaystyle  \Vert error \Vert \)},
ymajorgrids,
xtick={-6,-7,-8,-9,-10,-11},
xticklabels={$2^{-6}$,$2^{-7}$,$2^{-8}$,$2^{-9}$,$2^{-10}$,$2^{-11}$,},
ytick={-18,-16,-14,-12,-10,-8,-6},
yticklabels={$2^{-18}$,$2^{-16}$,$2^{-14}$,$2^{-12}$,$2^{-10}$,$2^{-8}$,$2^{-6}$},
ymin=-18, ymax=-5,
ytick style={color=black},
y post scale=1.0
]
\addplot [semithick, brown, mark=square*, mark size=3, mark options={solid}]
table {%
-6 -7.03492327720682
-7 -8.87893639303829
-8 -10.450046871219
-9 -11.75371440677
-10 -12.8866009819326
-11 -13.9419679739294
};
\addlegendentry{height}
\addplot [semithick, purple, mark=*, mark size=3, mark options={solid}]
table {%
-6 -6.30301696982099
-7 -8.22974973137303
-8 -9.81679448839414
-9 -10.9907586007867
-10 -11.9926524514412
-11 -12.9682485170229
};
\addlegendentry{velocity}
\addplot [thick, blue]
table {%
-9 -14
-9 -16
-11 -16
-9 -14
};
\node[left,blue] at (-9,-15) {2};
\node[below,blue] at (-10,-16) {2};

\nextgroupplot[
legend cell align={left},
legend style={fill opacity=0.8, draw opacity=1, text opacity=1, draw=lightgray204},
legend pos = south west,
scaled y ticks=manual:{}{\pgfmathparse{#1}},
tick align=outside,
tick pos=left,
title={$\alpha = 0$},
x dir=reverse,
x grid style={darkgray176},
xlabel={\(\displaystyle \Vert \Delta x \Vert \)},
xmajorgrids,
xmin=-11.25, xmax=-5.75,
xtick style={color=black},
y grid style={darkgray176},
ymajorgrids,
ytick distance=2,
ymin=-18, ymax=-5,
ytick style={color=black},
xtick={-6,-7,-8,-9,-10,-11},
xticklabels={$2^{-6}$,$2^{-7}$,$2^{-8}$,$2^{-9}$,$2^{-10}$,$2^{-11}$,},
ytick={-18,-16,-14,-12,-10,-8,-6},
yticklabels={$2^{-18}$,$2^{-16}$,$2^{-14}$,$2^{-12}$,$2^{-10}$,$2^{-8}$,$2^{-6}$},
y post scale=1.0,
]

\addplot [semithick, brown, mark=square*, mark size=3, mark options={solid}]
table {%
-6 -7.22053704046984
-7 -9.29141268782456
-8 -11.3126987277908
-9 -13.3188218070618
-10 -15.3193025513295
-11 -17.3200223393786
};
\addlegendentry{height}
\addplot [semithick, purple, mark=*, mark size=3, mark options={solid}]
table {%
-6 -6.1633893825145
-7 -8.18755572230195
-8 -10.1942610845455
-9 -12.1964686623078
-10 -14.1969677796315
-11 -16.1973031812642
};
\addlegendentry{velocity}
\addplot [thick, blue]
table {%
-9 -14
-9 -16
-10 -16
-9 -14
};
\node[left,blue] at (-9,-15) {2};
\node[below,blue] at (-9.5,-16) {1};
\end{groupplot}

\end{tikzpicture}